\documentclass[12pt]{article}
\usepackage{amsmath}
\usepackage{amsfonts}
\usepackage{amssymb}
\usepackage{amsthm}

\normalbaselines \oddsidemargin 0.5cm \evensidemargin 0.5cm
\usepackage{color}

\begin{document}
\pagestyle{empty}

\rightline{}

 \vskip 1cm
\begin{center}

{\Large {\textbf{New kinds of deformed Bessel functions}}}

\vspace{6mm}

{\bf \large   Mohammed Brahim Zahaf $^a$\footnote{E-mail address~:
m\_b\_zahaf@yahoo.fr} and Dominique Manchon $^{b}$\footnote{E-mail
address~: Dominique.Manchon@math.univ-bpclermont.fr}}
 \vspace{2mm}

{\small {\it $^{a}$ Laboratoire de Physique Quantique de
la Mati\`ere et Mod\'elisations Math\'ematiques (LPQ3M),
\\Universit\'e de Mascara, 29000-Mascara, Alg\'erie}}

\smallskip
{\small {\it $^{b}$Laboratoire de Math\'ematiques, CNRS-UMR 6620\\
Universit\'e Blaise Pascal\\
BP 80026, 63171 Aubi\`ere Cedex,\\
France }}

\end{center}
\vspace{4cm}
\begin{abstract}
{\small \noindent 
Using a  deformed
calculus based on the Dunkl operator, two new deformations of Bessel functions are proposed. Some properties i.e. generating function,  differential-difference equation, recursive relations, Poisson formula... are also given. Three more deformations are also outlined in the last section.
}
\end{abstract}
\pagestyle{plain} \setcounter{page}{1}
\section*{Introduction}
Dunkl operators on ${\mathbb R}^n$ have been introduced in \cite{Dunkl}: roughly speaking, they are partial derivative operators perturbed by reflexions. Although they are not differential operators in the usual sense, they mutually commute as classical partial derivative operators do. They have been intensively studied both from an algebraic \cite{DunklJeu, DunklOpdam, Torossian} and analytic \cite{trimeche1, Bensaid} point of view. E. M. Opdam proposed a family of deformed Bessel functions in this general context, which happens to coincide with the ordinary ones in the one-dimensional case \cite[Definition 6.9]{opdam}.\\

In this article, we exhibit several families of deformed Bessel functions in the one-dimensional case. We replace the derivation operator $d/dx$ by the corresponding Dunkl operator:
\begin{equation}
D_{\mu}f(x) = \frac{d}{dx}f(x)+\mu\frac{f(x)-f(-x)}{x}
\end{equation}
where $\mu$ is a real parameter, and consider the corresponding intertwining operator $V_\mu$. It is the unique linear operator such that $V_\mu\circ (d/dx)=D_\mu\circ V_\mu$ and $V_\mu(1)=1$ where $1$ is the constant function equal to one. Applying this operator to the generating function of the classical Bessel functions $J_n$, we obtain a first family  $J_n^\mu$ of deformed Bessel functions. They verify a differential-difference equation of order 3, as well as recursive formulas and a Poisson formula. They also verify a deformed version of the ``addition theorem'' \eqref{addition}, namely:
\begin{equation}
\tau _yJ_{n}^{\mu}(x)=\sum_{k=-\infty}^{\infty}J_{k}^{\mu}(x)J_{n-k}^{\mu}(y),
\end{equation}
where 
\begin{equation}
\tau _{y} f:=\sum_{n=0}^{\infty}\frac{y^n}{[n]_{\mu}!}D_{\mu}^n f,
\end{equation}
with $[n]_\mu=n+\mu\big(1-(-1)^n\big)$ and $[n]_\mu!=[n]_\mu [n-1]_\mu \cdots [1]_\mu$. Moreover, a connection formula relating the ${ J}_n^\mu$'s with their classical counterparts is easily given. Letting $\mu$ going to zero gives back the classical Bessel functions.  We also give a second family ${\mathcal J}_n^\mu$ of deformed Bessel functions by directly modifying the Poisson formula \eqref{poisson}. They also verify a differential-difference equation of order 3, but it seems there is no suitable addition theorem for them. Recursive relations are fulfilled, the generating function and the connection formula  are explicitly given. \\

We give three other families of deformed Bessel functions in the last section, together with the corresponding generating functions. They are obtained by directly deforming the coefficients of the Taylor series expansion of the $J_n$'s at the origin. Two among them verify recursive relations, but other properties of the classical Bessel functions do not seem to have their counterparts for these deformations.\\

\noindent
\textbf{Acknowledgements}: One of us (M. B. Z.) would like to thank Prof. B. Abdesselem and Dr. A. Yanallah for precious help and useful discussions. This work is partially granted by the project ANDRU/PNR/49/04/2011.
\section{Dunkl operator}
\noindent The Dunkl operator $D_{\mu}$ of index $\mu$ , $\mu \geq 0$, is defined on all smooth functions $f$ on $\mathbb{R}$ by 
\begin{equation}
D_{\mu}f(x) = \frac{d}{dx}f(x)+\mu\frac{f(x)-f(-x)}{x},\qquad
x\in\mathbb{R}.
\end{equation}
For
more general Dunkl operators see \cite{DunklOpdam}.
This operator has the following properties
\begin{eqnarray}
D^2_{\mu}f(x)&=&\frac{d^2}{dx^2}f(x)+\frac{2\mu}{x}\frac{d}{dx}f(x)-\frac{\mu}{x^2}(f(x)-f(-x)),\\
D_{\mu}(fg)(x) &=& f(x)D_{\mu}g(x)+g(-x)D_{\mu}f(x)+f'(x)(g(x)-g(-x)), \\
D_{\mu}x^n&=&\left[n\right]_{\mu}x^{n-1},\,\,n\in \mathbb{N}.
\end{eqnarray}
where
$\left[n\right]_{\mu}=n+\mu(1-(-1)^n)$. Obviously, $[2m]_{\mu}= 2m$, $[2m+1]_{\mu}=2m+1+2\mu$ for any integer $m$, and when $\mu\rightarrow 0$, $ [n]_{\mu}\rightarrow n$.
Let us define the deformed expenontial function by
\begin{equation}
E_{\mu}(x) =\sum_{n\geq 0}\frac{x^n}{\left[n\right]_{\mu}!},
\end{equation}
where $\left[n\right]_{\mu}!=\left[n\right]_{\mu}\left[n-1\right]_{\mu}...\left[1\right]_{\mu}$,\qquad
$\left[0\right]_{\mu}!\equiv 1$. Then we have 
\begin{equation}
D_{\mu}E_{\mu}(\lambda x)=\lambda E_{\mu}(\lambda x),\,\,\lambda\in \mathbb{C}.
\end{equation}
Recall the definition of the Pochhammer symbol:
\begin{equation}
(a)_k=\frac{\Gamma(a +k)}{\Gamma (a)}=a(a+1)\cdots(a+k-1).
\end{equation}
Let us remark, using the following expressions \cite{rosenblum}:
\begin{eqnarray}
\left[2m\right]_{\mu}!&= &\frac{2^{2m}m!\Gamma(m+\mu+\frac{1}{2})}{\Gamma(\mu+\frac{1}{2})}=(2m)!\frac{(\mu+\frac{1 }{2})_{m}}{(\frac{1 }{2})_{m}},\nonumber\\
\left[2m+1\right]_{\mu}!&=& \frac{2^{2m+1}m!\Gamma(m+\mu+\frac{3}{2})}{\Gamma(\mu+\frac{1}{2})}=(2m+1)!\frac{(\mu+\frac{1 }{2})_{m+1}}{(\frac{1 }{2})_{m+1}},
\end{eqnarray}
that we can easily write  $E_{\mu}(x)$ under the form
\begin{eqnarray}
E_{\mu}(x) &=&j_{\mu-\frac{1}{2}}(ix)+\frac{x}{2\mu+1}j_{\mu+\frac{1}{2}}(ix)\nonumber\\
&=&e^{x}\,   _1F_1\left(
\begin{array}{l}
\mu\\
2\mu+1
\end{array};-2x\right),
\end{eqnarray}
where $j_{\alpha}$ is the normalized spherical Bessel function defined for $\alpha \geq -\frac{1}{2}$, by
\begin{equation}
j_{\alpha}(x)=\sum_{k\geq 0}\frac{(-1)^k}{k!(\alpha +1)_k}\left(\frac{x}{2}\right)^{2k}. \nonumber
\end{equation}
For our purpose let us recall the following important theorem (\cite{Dunkl}, \cite{DunklJeu}): \textsl {There exists a unique linear isomorphism $V_{\mu}$ (called Dunkl intertwining operator) from the set of polynomials $\mathcal{P}_n$ on $\mathbb{R}$ of degree less or equal than $n$ onto itself such that:}
\begin{equation}\label{dun}
 V_{\mu}(1)=1,\quad and\,\,\,\, D_{\mu}V_{\mu}=V_{\mu}\frac{d}{dx}.
\end{equation}
The operator $V_{\mu}$ has been extended by K. Trim\`eche to an
isomorphism from $\mathcal{C}^{\infty}(\mathbb{R})$ onto itself satisfying
the relations in \eqref{dun} (see \cite{trimeche1}).
It  possesses the following integral representation:
\begin{eqnarray}
\forall x\in \mathbb{R},\,\,\,V_{\mu}(f(x))=\frac{1}{\beta(\frac{1}{2},\mu)}\int_{-1}^{1}f(xt)(1-t)^{\mu-1}(1+t)^{\mu}dt ,\,\,\,\,f\in\mathcal{C}^{\infty}(\mathbb{R}).
\end{eqnarray}
We have
\begin{eqnarray}\label{vn}
V_{\mu}(x^n)&=&\frac{(\frac{1}{2})_{\lfloor\frac{n+1}{2}\rfloor}}{(\mu+\frac{1}{2})_{\lfloor\frac{n+1}{2}\rfloor}}x^n=
\frac{n!}{[n]_{\mu}!}x^n,
\end{eqnarray}
where $\lfloor \alpha\rfloor$ stands for  integer  part of the real number $\alpha$, 
and
\begin{equation}\label{Evn}
E_{\mu}(x)=V_{\mu}(e^x).
\end{equation}
The generalized translation operator $\tau _y$, $y\in \mathbb{R}$ is defined by
\begin{equation}
\tau _{y} f:=E_{\mu}(yD_{\mu})f=\sum_{n=0}^{\infty}\frac{y^n}{[n]_{\mu}!}D_{\mu}^n f.
\end{equation}
for all entire functions $f$ on $\mathbb{C}$ for which the series converges pointwise. It possesses the following properties:
 \begin{eqnarray}
\tau _{y}x^n=\sum_{k=0}^{+\infty}\frac{y^k}{[k]_{\mu}!}D_{\mu}^kx^n=\sum_{k=0}^{+\infty}\left(
\begin{array}{c}
n\\
k
\end{array}
\right)_{\mu}x^ky^{n-k},
\end{eqnarray}
where $\left(
\begin{array}{c}
n\\
k
\end{array}
\right)_{\mu}=\frac{[n]_{\mu}!}{[k]_{\mu}![n-k]_{\mu}!}$. We moreover have:
\begin{equation}
\tau _yE_{\mu}(\lambda x)=E_{\mu}(\lambda x)E_{\mu}(\lambda y), \,\, \lambda \in \mathbb{C}.
\end{equation}
\section{Background on the classical Bessel functions}
Let $n\in{\mathbb Z}$ be any integer. The classical Bessel function of order $n$ is given by:
\begin{eqnarray}
J_n(x)&=&\sum_{k\geq 0}\frac{(-1)^k}{k!(k+n)!}\left(\frac{x}{2}\right)^{2k+n}\nonumber\\
&=& \frac{1}{n!}\left(\frac{x}{2}\right)^n\medskip  \medskip _{0}F_{1}
\left(
\begin{array}{l}
-\\
n+1
\end{array}
;-\frac{x^2}{4}\right)
.
\end{eqnarray}
In view of the relation:
\begin{equation}\label{parity}
J_{-n}(x)=(-1)^nJ_n(x),
\end{equation}
we shall consider only $J_n$ for nonnegative $n$. The order $n$ Bessel function is the solution of the following second-order linear differential equation:
\begin{equation}\label{bessel-equadif}
\left(x^2 \dfrac{d^2}{dx^2}+x\dfrac{d}{dx}+(x^2-n^2)\right)y(x)=0.
\end{equation}
with boundary conditions $J_n(0)=\delta_0^n$ and $\dot J_n(0)=\frac {1}{2}\delta_n^1$.
The classical Bessel functions can be gathered into the generating function:
\begin{eqnarray}
G(x,t)=\exp\left(\frac{x}{2}(t-\frac{1}{t})\right)&=&\sum_{n=-\infty}^{+\infty}J_n(x)t^n.
\end{eqnarray}
The following recursive relations are satisfied:
\begin{eqnarray}
2J'_n(x)&=&J_{n-1}(x)-J_{n+1}(x),\label{bessel1}\\
nJ_n(x)&=&xJ_{n-1}(x)-xJ'_{n}(x),\label{bessel2}\\
nJ_n(x)&=&xJ_{n+1}(x)+x J'_{n}(x),\label{bessel3}
\end{eqnarray}
as well as
\begin{equation}\label{bessel4}
J_{n-1}(x)+J_{n+1}(x)=\frac{2n}{x}J_n(x),
\end{equation}
which can be obtained by adding \eqref{bessel2} and \eqref{bessel3}. The \textsl{Poisson formula} is given by:
\begin{equation}\label{poisson}
J_{n}(x)=\frac{(\frac{x}{2})^n}{\Gamma(\frac{1}{2})\Gamma(n+\frac{1}{2})}\int_{-1}^{1}(1-s^2)^{n-\frac{1}{2}}\cos(sx)ds,
\end{equation}
The following \textsl{addition theorem} holds:
\begin{equation}\label{addition}
J_n(x+y)=\sum_{k=-\infty}^{+\infty}J_k(x)J_{n-k}(y),
\end{equation}
as one can easily see by using the generating function $G(x+y, t)$. For more details on classical Bessel functions, see for example  \cite{AS, lebedev, Watson}.
\section{First deformation of the Bessel function}
We define the deformed Bessel function by
\begin{equation}
J_{n}^{\mu}(x):=V_{\mu}\big(J_n(x)\big),
\end{equation}
where $J_n$ is the classical Bessel function. In virtue of  \eqref{vn} we can write
\begin{equation}\label{defbessel}
J_{n}^{\mu}(x)=\sum_{k \geq 0}\frac{(-1)^k (2k+n)!}{k!(k+n)!\left[2k+n\right]_{\mu}!}\left(\frac{x}{2}\right)^{2k+n}.
\end{equation}
Using the ratio test, we can verify that this series converges in the whole
complex plane, and hence represents an entire function of $x$. Letting $V_{\mu}$ operate on both sides of \eqref{parity} we deduce that:
\begin{equation}
J_{-n}^{\mu}(x)=(-1)^{n}J_{n}^{\mu}(x),\,\,n=1,2,...
\end{equation}

The deformed Bessel function of
order $n$ is a solution of the following differential-difference equation based on the deformed derivative (Dunkl operator):
\begin{equation}\label{eqdiff-un}
\Big((xD_{\mu}-[n]_{\mu})(xD_{\mu}-[-n]_{\mu})(xD_{\mu}+\beta _{n}+2\mu-1)+x^2(xD_{\mu}+\beta _{n}+1)\Big)y(x)=0,\nonumber\\
\end{equation}
where 
\begin{eqnarray}\label{betan}
\beta_{n}=2\lfloor\frac{n+1}{2}\rfloor-[n]_{\mu}=\left\{
\begin{array}{ll}
0& {\text {if}}\,\, n\in 2\mathbb{N}\\
1-2\mu & {\text {if}}\,\,n\in 2\mathbb{N}+1
\end{array}
\right..
\end{eqnarray}

In fact,
\begin{eqnarray}
\begin{array}{l}
\Big((xD_{\mu}-[n]_{\mu})(xD_{\mu}-[-n]_{\mu})(xD_{\mu}+\beta _{n}+2\mu-1)\Big)J_{n}^{\mu}(x)=\nonumber\\
=\sum_{k \geq 0}\frac{(-1)^k (2k+n)!}{k!(k+n)!\left[2k+n\right]_{\mu}!}\Big((2k+[n]_{\mu}-[n]_{\mu})(2k+[n]_{\mu}-[-n]_{\mu})(2k+[n]_{\mu}+\beta_{n}+2\mu-1)\Big)\left(\frac{x}{2}\right)^{2k+n}\nonumber\\
=\sum_{k \geq 0}\frac{(-1)^k (2k+n)!}{k!(k+n)!\left[2k+n\right]_{\mu}!}\Big((2k)(2k+2n)(2k+2{\lfloor\frac{n+1}{2}\rfloor}+2\mu-1)\Big)\left(\frac{x}{2}\right)^{2k+n}\nonumber\\
=\sum_{k \geq 0}\frac{(-1)^k (2k+n)!}{(k-1)!(k+n-1)!\left[2k+n\right]_{\mu}!}\Big(4(2k+2{\lfloor\frac{n+1}{2}\rfloor}+2\mu-1)\Big)\left(\frac{x}{2}\right)^{2k+n}\nonumber\\
=-x^2\sum_{k \geq 0}\frac{(-1)^k (2k+n+2)!}{k!(k+n)!\left[2k+n+2\right]_{\mu}!}\Big((2k+2+2{\lfloor\frac{n+1}{2}\rfloor}+2\mu-1)\Big)\left(\frac{x}{2}\right)^{2k+n}\nonumber\\
=-x^2\sum_{k \geq 0}\frac{(-1)^k (2k+n)!}{k!(k+n)!\left[2k+n\right]_{\mu}!}\left(\frac{(2k+n+2)(2k+n+1)}{\left[2k+n+2\right]_{\mu}\left[2k+n+1\right]_{\mu}}(2k+2
\lfloor\frac{n+1}{2}\rfloor+2\mu+1)\right)\left(\frac{x}{2}\right)^{2k+n}\nonumber\\
=-x^2\sum_{k \geq 0}\frac{(-1)^k (2k+n)!}{k!(k+n)!\left[2k+n\right]_{\mu}!}\Big((2k+2{\lfloor\frac{n+1}{2}\rfloor}+1)\Big)\left(\frac{x}{2}\right)^{2k+n}\nonumber\\
=-\Big(x^2(xD_{\mu}+\beta _{n}+1)\Big)J_{n}^{\mu}(x).
\end{array}
\end{eqnarray}
It is easy to see, when $\mu=0$, that the third order deformed differential equation \eqref{eqdiff-un} reduces to the second order Bessel differential equation \eqref{bessel-equadif}.
The  generating function of the deformed Bessel function is given by:
\begin{equation}\label{gen}
G^{\mu}(x,t)=E_{\mu}\left(\frac{x}{2}(t-\frac{1}{t})\right)=\sum_{n=-\infty}^{+\infty}J_n^{\mu}(x)t^n.
\end{equation}
This is obtained by applying the intertwining operator $V_{\mu}$ to  the generating function $G(x,t)$ for the  classical Bessel function, with respect to the variable $x$, and using the relation \eqref{Evn}. For $t=1$, we obtain the  following relation:
\begin{equation}
\sum_{n=-\infty}^{+\infty}J_n^{\mu}(x)=1,
\end{equation}
which can be also writen as:
\begin{equation}
J_0^{\mu}(x)+2\sum_{n=1}^{+\infty}J_{2n}^{\mu}(x)=1.
\end{equation}
If we take $t=e^{i\theta}$ in \eqref{gen} we obtain:
\begin{equation}\label{gen2}
E_{\mu}(ix\sin\theta)=\sum_{n=-\infty}^{+\infty}J_n^{\mu}(x)e^{in\theta}.
\end{equation}
This implies that:
\begin{equation}
J_n^{\mu}(x)=\frac{1}{2\pi}\int_{-\pi}^{\pi}E_{\mu}(ix\sin\theta)e^{-in\theta}d\theta.
\end{equation}
The deformed Bessel function possesses the  following recurrence relations:
\begin{eqnarray}
2D_{\mu}J_n^{\mu}(x)&=&J_{n-1}^{\mu}(x)-J_{n+1}^{\mu}(x),\label{dbessel1}\\
nD_{\mu}J_n^{\mu}(x)&=&\frac{d}{dx}\left(xJ_{n-1}^{\mu}(x)-xD_{\mu}J_{n}^{\mu}(x)\right),\label{dbessel2}\\
nD_{\mu}J_n^{\mu}(x)&=&\frac{d}{dx}\left(xJ_{n+1}^{\mu}(x)+xD_{\mu}J_{n}^{\mu}(x)\right).\label{dbessel3}
\end{eqnarray}
 Summing up \eqref{dbessel2} and \eqref{dbessel3} we obtain:
\begin{equation}\label{dbessel4}
2nD_\mu J_n^\mu(x)=\frac d{dx}\big(xJ_{n-1}^\mu(x)+xJ_{n+1}^\mu(x)\big).
\end{equation}
The three last equations \eqref{dbessel2},  \eqref{dbessel3}  and \eqref{dbessel4} are not as simple as their classical counterparts \eqref{bessel2},  \eqref{bessel3}  and \eqref{bessel4}, due to the fact that the deformed derivative $D_\mu$ differs from the ordinary derivative $d/dx$. The first relation is obtained by applying the intertwining operator to the left and the right hand of (\ref{bessel1}). For the second one, we have:
\begin{eqnarray}
J_{n-1}^{\mu}(x)-D_{\mu}J_{n}^{\mu}(x)&=&V_{\mu}(J_{n-1}(x))-D_{\mu}V_{\mu}(J_{n}(x))\nonumber\\
&=&V_{\mu}(J_{n-1}(x)-J'_{n}(x))\nonumber\\
&=&V_{\mu}\left(\frac{n}{x}J_{n}(x)\right),\nonumber
\end{eqnarray}
where we  have used  the fact that $D_{\mu}V_{\mu}=V_{\mu}\frac{d}{dx}$ and the relation \eqref{bessel2}, therefore
\begin{eqnarray}
\frac{d}{dx}\left(xJ_{n-1}^{\mu}(x)-xD_{\mu}J_{n}^{\mu}(x)\right)&=&\frac{d}{dx}\left(x\left(V_{\mu}\left(\frac{n}{x}J_{n}(x)\right)\right)\right)\nonumber\\
&=&n\frac{d}{dx}\left(\frac{x}{\beta(\frac{1}{2},\mu)}\int_{-1}^{1}\frac{J_n(xt)}{xt}(1-t)^{\mu-1}(1+t)^{\mu}dt \right)\nonumber\\
&=&n\frac{1}{\beta(\frac{1}{2},\mu)}\frac{d}{dx}\left(\int_{-1}^{1}\frac{J_n(xt)}{t}(1-t)^{\mu-1}(1+t)^{\mu}dt \right)\nonumber\\
&=&n\frac{1}{\beta(\frac{1}{2},\mu)}\int_{-1}^{1}J'_n(xt)(1-t)^{\mu-1}(1+t)^{\mu}dt \nonumber\\
&=&nV_{\mu}(J'_n(x))=nD_{\mu}J_n^\mu(x).\nonumber
\end{eqnarray}
The third one follows similarly.\\

A formula involving deformed Bessel functions with different {\em superscripts} is
\begin{equation}
(xD_{\mu}+\beta _{n}+2\mu-1)J_{n}^{\mu}(x)=(2\mu-1)J_{n}^{\mu-1}(x),
\end{equation}
where $\beta_{n}$ is given in \eqref{betan}. To prove this equality, let first remark that
\begin{eqnarray}\nonumber
\beta_{n}+2\mu-1= \left\{
\begin{array}{ll}
2\mu-1& n\in 2\mathbb{N}\\
0& n\in 2\mathbb{N}+1
\end{array}
\right.
\end{eqnarray}
and
\begin{equation}
\frac{2(\mu-1)}{\beta(\frac{1}{2},\mu)}=\frac{2\mu-1}{\beta(\frac{1}{2},\mu-1)}.\nonumber
\end{equation}
In virtue of \eqref{dun} we have
\begin{equation}
xD_{\mu}J_{n}^{\mu}(x)=xV_{\mu}(J'_n(x))=\frac{1}{\beta(\frac{1}{2},\mu)}\int_{-1}^{1}xJ'_n(xt)(1-t)^{\mu-1}(1+t)^{\mu}dt,\nonumber
\end{equation}
by using an integration by parts, we have

\begin{eqnarray}\label{int}
xD_{\mu}J_{n}^{\mu}(x)&=&\frac{1-2\mu}{\beta(\frac{1}{2},\mu)}\int_{-1}^{1}J_n(xt)(1-t^2)^{\mu-1}dt+\frac{2(\mu-1)}{\beta(\frac{1}{2},\mu)}\int_{-1}^{1}J_n(xt)(1-t)^{\mu-2}(1+t)^{\mu-1}dt\nonumber\\
&=&\frac{1-2\mu}{\beta(\frac{1}{2},\mu)}\int_{-1}^{1}J_n(xt)(1-t^2)^{\mu-1}dt+(2\mu-1)J_n^{\mu-1}(x).
\end{eqnarray}
If $n\in  2\mathbb{N}$, the function $J_n(x)$ is even and 
\begin{eqnarray}
\frac{1-2\mu}{\beta(\frac{1}{2},\mu)}\int_{-1}^{1}J_n(xt)(1-t^2)^{\mu-1}dt=\frac{1-2\mu}{\beta(\frac{1}{2},\mu)}\int_{-1}^{1}J_n(xt)(1+t)(1-t^2)^{\mu-1}dt=-(2\mu -1)J_{n}^{\mu}(x),\nonumber
\end{eqnarray}
if $n\in  2\mathbb{N}+1$, the function $J_n(x)$ is odd and the integrale in \eqref{int} is equal to $0$, and therefore: 
\begin{equation}
xD_{\mu}J_{n}^{\mu}(x)=-(\beta _{n}+2\mu-1)J_n^{\mu}(x)+(2\mu-1)J_n^{\mu-1}(x).\nonumber
\end{equation}

Using the fact:
\begin{eqnarray}
\frac{(2k+n)!}{[2k+n]_{\mu}!}&=&\frac{(\frac{1}{2})_{k+\lfloor\frac{n+1}{2}\rfloor}}{(\mu+\frac{1}{2})_{k+\lfloor\frac{n+1}{2}\rfloor}}=\frac{(\frac{1}{2})_{\lfloor\frac{n+1}{2}\rfloor}}{(\mu+\frac{1}{2})_{\lfloor\frac{n+1}{2}\rfloor}}\frac{(\lfloor\frac{n+1}{2}\rfloor+\frac{1}{2})_{k}}{(\lfloor\frac{n+1}{2}\rfloor+\mu+\frac{1}{2})_{k}}\nonumber\\&=&\frac{n!}{[n]_{\mu}!}\frac{(\lfloor\frac{n+1}{2}\rfloor+\frac{1}{2})_{k}}{(\lfloor\frac{n+1}{2}\rfloor+\mu+\frac{1}{2})_{k}},
\end{eqnarray}
the deformed Bessel function can be expressed in terms of the generalized hypergeometric function $\medskip _{1}F_{2}$  as follows:
\begin{eqnarray}
J_{n}^{\mu}(x)&=&\frac{n!}{[n]_{\mu}!}\sum_{k \geq 0}\frac{(-1)^k }{k!(k+n)!}\frac{(\lfloor\frac{n+1}{2}\rfloor+\frac{1}{2})_{k}}{(\lfloor\frac{n+1}{2}\rfloor+\mu+\frac{1}{2})_{k}}\left(\frac{x}{2}\right)^{2k+n}\nonumber\\
&=&\frac{1}{[n]_{\mu}!}\left(\frac{x}{2}\right)^n\sum_{k \geq 0}\frac{(\lfloor\frac{n+1}{2}\rfloor+\frac{1}{2})_{k}}{k!(n+1)_k(\lfloor\frac{n+1}{2}\rfloor+\mu+\frac{1}{2})_{k}}\left(-\frac{x^2}{4}\right)^{k}\nonumber\\
&=&\frac{1}{[n]_{\mu}!}\left(\frac{x}{2}\right)^n\medskip  \medskip _{1}F_{2}
\left(
\begin{array}{ll}
\lfloor\frac{n+1}{2}\rfloor+\frac{1}{2}\\
n+1&\lfloor\frac{n+1}{2}\rfloor+\mu+\frac{1}{2}
\end{array}
;-\frac{x^2}{4}\right).
\end{eqnarray}
By applying the  intertwining operator $V_{\mu}$ to the left and right hand of the  classical  Bessel function Poisson formula \eqref{poisson}, we obtain the $\mu$-version of Poisson formula
\begin{eqnarray}
J_{n}^{\mu}(x)&=&\frac{n!}{[n]_{\mu}!}\frac{(\frac{x}{2})^n}{\Gamma(\frac{1}{2})\Gamma(n+\frac{1}{2})}\int_{-1}^{1}(1-s^2)^{n-\frac{1}{2}}\medskip  \medskip _{1}F_{2}
\left(
\begin{array}{ll}
\lfloor\frac{n+1}{2}\rfloor+\frac{1}{2}\\
\frac{1}{2}&\lfloor\frac{n+1}{2}\rfloor+\mu+\frac{1}{2}
\end{array}
;-\frac{(sx)^2}{4}\right)ds.\nonumber\\
\end{eqnarray} 
The deformed Bessel function $J_{n}^{\mu}(x)$ verifies the following analogue of the addition theorem \eqref{addition}:
\begin{equation}\label{daddition}
\tau _yJ_{n}^{\mu}(x)=\sum_{k=-\infty}^{\infty}J_{k}^{\mu}(x)J_{n-k}^{\mu}(y).
\end{equation}
This follows immediately from the fact that:
\begin{equation}
\tau _yE_{\mu}\left(\frac{x}{2}(t-\frac{1}{t})\right)=E_{\mu}\left(\frac{x}{2}(t-\frac{1}{t})\right)E_{\mu}\left(\frac{y}{2}(t-\frac{1}{t})\right),\nonumber
\end{equation}
or equivalently:
\begin{equation}
\sum_{n=-\infty}^{+\infty}\tau _yJ_n^{\mu}(x)t^n =\sum_{n=-\infty}^{+\infty}J_n^{\mu}(x)t^n\sum_{n=-\infty}^{+\infty}J_n^{\mu}(y)t^n.\nonumber
\end{equation}
The result in \eqref{daddition} is obtained by equating the coefficients of $t^n$.\\

A connection formula between the deformed Bessel function ${J}_{n}^{\mu}$ and the classical one is given by:
\begin{equation}\label{connection1}
{J}_{n}^{\mu}(x)=\frac{n!}{[n]_{\mu}!}\sum_{k\geq 0}\frac{(\mu)_k}{k!(\lfloor\frac{n+1}{2}\rfloor+\mu+\frac{1}{2})_k}\left(\frac{x}{2}\right)^k{J}_{n+k}(x).
\end{equation}
To prove this formula we use the fact
\begin{eqnarray}
\frac{(a)_k}{(b)_k}=\medskip _{2}F_1
\left(
\begin{array}{cl}
-k&b-a\\
&b
\end{array}
;1\right)=\sum_{k'=0}^{k}\frac{(-k)_{k'}(b-a)_{k'}}{k'!(b)_{k'}}=\sum_{k'=0}^{k}\frac{(-1)^{k'}k!(b-a)_{k'}}{k'!(k-k')!(b)_{k'}},
\end{eqnarray}
which leads to
\begin{eqnarray}\label{gauss}
\frac{(2k+n)!}{[2k+n]_{\mu}!}&=&\frac{n!}{[n]_{\mu}!}\sum_{k'=0}^{k}\frac{(-1)^{k'}k!(\mu)_{k'}}{k'!(k-k')!(\lfloor\frac{n+1}{2}\rfloor+\mu+\frac{1}{2})_{k'}},
\end{eqnarray}
and therefore the formula  \eqref{connection1} is obtained after substituting the right hand in \eqref{gauss} in the series \eqref{defbessel}.

The deformed Bessel functions of order $n$ have the following asymptotic form, as $|x|\rightarrow \infty$ and  $|\arg(x)|\le\frac{\pi}{2}-\varepsilon$ for some $\varepsilon>0$:
\begin{eqnarray}
J_n^{\mu}(x)&\approx&\frac{\Gamma(\mu+\frac{1}{2})}{\pi}\left(\frac{x}{2}\right)^{-\mu-\frac{1}{2}}\cos\left(x-\frac{\pi}{2}(\frac{1}{2}+n+\mu)\right)\nonumber\\
&+&\frac{\Gamma(\lfloor\frac{n+1}{2}\rfloor+\frac{1}{2})\Gamma(\mu+\frac{1}{2})}{\Gamma(n+\frac{1}{2}-\lfloor\frac{n+1}{2}\rfloor)\Gamma(\mu)\Gamma(\frac{1}{2})}\left(\frac{x}{2}\right)^{n-2\lfloor\frac{n+1}{2}\rfloor-1},
\end{eqnarray}
which is obtained by using the asymptotic formula of the generalized hypergeometric function $\medskip _{1}F_{2}$ :
\begin{eqnarray}
 \medskip _{1}F_{2}
\left(
\begin{array}{ll}
\alpha\\
n+1&\beta
\end{array}
;-x^2\right)\approx \frac{\Gamma(\beta)n!}{\Gamma(\frac{1}{2})\Gamma(\alpha)}x^{-\frac{1}{2}-n+\alpha-\beta}\cos\left(2x-\frac{\pi}{2}\left(\frac{1}{2}+n+\beta-\alpha\right)\right)\nonumber\\+\frac{\Gamma(\beta)n!}{\Gamma(\beta-\alpha)\Gamma(n+1-\alpha)}x^{-2\alpha},
\end{eqnarray}
as $|x|\rightarrow \infty$ and $|\arg(x)|\le\frac{\pi}{2}-\varepsilon$ for some $\varepsilon>0$,
which is a special case of a general formula given by Luke \cite[ p. 203, Eq. (4)]{luke}.

\section{A second deformation of the Bessel function}
We define a second deformation of Bessel function by using an other deformation of the Poisson formula:
\begin{equation}
\mathcal{J}_n^{\mu}(x):=\frac{(\frac{x}{2})^n}{\Gamma(\frac{1}{2})\Gamma(n+\frac{1}{2})}\int_{-1}^{1}E_{\mu}(isx)(1-s^2)^{n-\frac{1}{2}}ds.
\end{equation}
Using the series representation of $E_{\mu}(isx)$ we can write the function $\mathcal{J}_n^{\mu}(x)$ as:
\begin{equation}
\mathcal{J}_n^{\mu}(x)=\sum_{k \geq 0}\frac{(-1)^k (2k)!}{k!(k+n)!\left[2k\right]_{\mu}!}\left(\frac{x}{2}\right)^{2k+n}.
\end{equation}
We remark that $\mathcal{J}_0^{\mu}(x)={J}_0^{\mu}(x)$. The deformed Bessel function $\mathcal{J}_n^{\mu}(x)$ is a solution of the deformed differential equation: 
\begin{eqnarray}\label{eqdiff}
\Big((xD_{\mu}-[n]_{\mu})(xD_{\mu}-[-n]_{\mu})(xD_{\mu}-[n]_{\mu}+2\mu-1)+x^2(xD_{\mu}-[n]_{\mu}+1)\Big)y(x)=0.\nonumber\\
\end{eqnarray}
In fact,
\begin{eqnarray}
\begin{array}{l}
\Big((xD_{\mu}-[n]_{\mu})(xD_{\mu}-[-n]_{\mu})(xD_{\mu}-[n]_{\mu}+2\mu-1)\Big)\mathcal{J}_{n}^{\mu}(x)=\nonumber\\
=\sum_{k \geq 0}\frac{(-1)^k (2k)!}{k!(k+n)!\left[2k\right]_{\mu}!}\Big((2k+[n]_{\mu}-[n]_{\mu})(2k+[n]_{\mu}-[-n]_{\mu})(2k+[n]_{\mu}-[n]_{\mu}+2\mu-1)\Big)\left(\frac{x}{2}\right)^{2k+n}\nonumber\\
=\sum_{k \geq 0}\frac{(-1)^k (2k)!}{k!(k+n)!\left[2k\right]_{\mu}!}\Big((2k)(2k+2n)(2k+2\mu-1)\Big)\left(\frac{x}{2}\right)^{2k+n}\nonumber\\
=\sum_{k \geq 0}\frac{(-1)^k (2k)!}{(k-1)!(k+n-1)!\left[2k\right]_{\mu}!}\Big(4(2k+2\mu-1)\Big)\left(\frac{x}{2}\right)^{2k+n}\nonumber\\
=-x^2\sum_{k \geq 0}\frac{(-1)^k (2k+2)!}{k!(k+n)!\left[2k+2\right]_{\mu}!}\Big((2k+2+2\mu-1)\Big)\left(\frac{x}{2}\right)^{2k+n}\nonumber\\
=-x^2\sum_{k \geq 0}\frac{(-1)^k (2k)!}{k!(k+n)!\left[2k\right]_{\mu}!}\Big(\frac{ (2k+2)(2k+1)}{\left[2k+2\right]_{\mu}\left[2k+1\right]_{\mu}}(2k+2\mu+1)\Big)\left(\frac{x}{2}\right)^{2k+n}\nonumber\\
=-x^2\sum_{k \geq 0}\frac{(-1)^k (2k)!}{k!(k+n)!\left[2k\right]_{\mu}!}(2k+1)\left(\frac{x}{2}\right)^{2k+n}\nonumber\\
=-x^2(xD_{\mu}-[n]_{\mu}+1)\mathcal{J}_{n}^{\mu}(x).
\end{array}
\end{eqnarray}
The  generating function of the deformed Bessel functions $\mathcal {J}_n^{\mu}(x)$ for $n\ge 0$ is given by:
\begin{eqnarray}\label{gen2}
\mathcal {G}^{\mu}(x,t)=e^{\frac{xt}{2}}\medskip  \medskip _{1}F_{1}
\left(
\begin{array}{l}
\frac{1}{2}\\
\mu+\frac{1}{2}
\end{array}
;-\frac{x}{2t}\right)&=&\mathcal {J}_0^{\mu}(x)+\sum_{n=1}^{+\infty}\mathcal {J}_n^{\mu}(x)\big(t^n+(-1)^nt^{-n}\big).
\end{eqnarray}
In fact, let us take
\begin{eqnarray}
e^{\frac{xt}{2}}\medskip  \medskip _{1}F_{1}
\left(
\begin{array}{l}
\frac{1}{2}\\
\mu+\frac{1}{2}
\end{array}
;-\frac{x}{2t}\right)&=&\sum_{n=-\infty}^{+\infty}c_n(x)t^n.
\end{eqnarray}
To calculate the coefficients $c_n(x)$, we multiply the power series
\begin{equation}
e^{\frac{xt}{2}}=\sum_{n\geq 0}\frac{(\frac{xt}{2})^{n}}{n!},
\end{equation}
\begin{eqnarray}
\medskip  \medskip _{1}F_{1}
\left(
\begin{array}{l}
\frac{1}{2}\\
\mu+\frac{1}{2}
\end{array}
;-\frac{x}{2t}\right)&=&\sum_{n\geq 0}\frac{(\frac{1}{2})_n}{n!(\mu+\frac{1}{2})_n}\left(-\frac{x}{2t}\right)^n,
\end{eqnarray}
and then combine terms containing identical powers of $t$. As a result, we
obtain
\begin{eqnarray}
c_n&=&\mathcal {J}_n^{\mu}(x),\,\,\,\,\,\,n=0,1,2,...,\nonumber\\
c_n&=&(-1)^n\mathcal {J}_{-n}^{\mu}(x),\,\,\,\,\,\,n=-1,-2,...,
\end{eqnarray}
which implies \eqref{gen2}. The deformed Bessel function $\mathcal{J}_{n}^{\mu}(x)$ possesses the  following recursive relations
\begin{equation}
(xD_{\mu}-[-n]_{\mu} )\mathcal{J}_{n}^{\mu}(x)=x\mathcal{J}_{n-1}^{\mu}(x),
\end{equation}
\begin{equation}\label{rec2}
\left(xD_{\mu}-[n]_{\mu}-1\right)\left(x\mathcal{J}_{n+1}^{\mu}(x)+x\mathcal{J}_{n-1}^{\mu}(x)-2(n-\mu) \mathcal{J}_{n}^{\mu}(x)\right)=-2\mu\mathcal{J}_{n}^{\mu}(x).
\end{equation}
Indeed, for the first one we have:
\begin{eqnarray}
(xD_{\mu}-[-n]_{\mu})\mathcal{J}_{n}^{\mu}(x)&=&\sum_{k \geq 0}\frac{(-1)^k (2k)!(2k+2n)}{k!(k+n)!\left[2k\right]_{\mu}!}\left(\frac{x}{2}\right)^{2k+n}\nonumber\\
&=&x\sum_{k \geq 0}\frac{(-1)^k (2k)!}{k!(k+n-1)!\left[2k\right]_{\mu}!}\left(\frac{x}{2}\right)^{2k+n-1}\nonumber\\
&=&x\mathcal{J}_{n-1}^{\mu}(x).\nonumber
\end{eqnarray}
For the second one, we have:
\begin{eqnarray}
&&x\mathcal{J}_{n+1}^{\mu}(x)+x\mathcal{J}_{n-1}^{\mu}(x)-2(n-\mu) \mathcal{J}_{n}^{\mu}(x)=\nonumber\\&=&\sum_{k \geq 0}\frac{(-1)^k (2k)!}{k!(k+n)!\left[2k\right]_{\mu}!}2\left(-\frac{k[2k-1]_{\mu}[2k]_{\mu}}{(2k-1)(2k)}+n+k-n+\mu \right)\left(\frac{x}{2}\right)^{2k+n}\nonumber
\\&=&\sum_{k \geq 0}\frac{(-1)^k (2k)!}{k!(k+n)!\left[2k\right]_{\mu}!}\left(-\frac{2\mu}{2k-1}\right)\left(\frac{x}{2}\right)^{2k+n},\nonumber
\end{eqnarray}
and after applying the operator $(xD_{\mu}-[n]_{\mu}-1)$, we retrieve the relation \eqref{rec2}.\\

Using the fact that:
\begin{eqnarray}
_1F_1\left(
\begin{array}{c}
a\\
b
\end{array};x\right)=e^{x}\,_1F_1\left(
\begin{array}{c}
b-a\\
b
\end{array};-x\right)
\end{eqnarray}

we can write:
\begin{eqnarray}
\mathcal {G}^{\mu}(x,t)={G}(x,t)_1F_1\left(
\begin{array}{c}
\mu\\
\mu+\frac{1}{2}
\end{array};\frac{x}{2t}\right)
\end{eqnarray}
which gives a connection formula between the deformed Bessel function and the classical one:
\begin{equation}\label{connection}
\mathcal{J}_{n}^{\mu}(x)=\sum_{k\geq 0}\frac{(\mu)_k}{k!(\mu+\frac{1}{2})_k}\left(\frac{x}{2}\right)^k{J}_{n+k}(x).
\end{equation}
By using this last relation and the recursive relations of the classical Bessel function \eqref{bessel1}-\eqref{bessel4} we can prove the following formulae involving deformed Bessel functions $\mathcal{J}_{n}^{\mu}(x)$ with different {\em superscripts} :
\begin{equation}
(xD_{\mu}-[n]_{\mu}+2\mu-1)\mathcal{J}_{n}^{\mu}(x)=(2\mu-1)\mathcal{J}_{n}^{\mu-1}(x),
\end{equation}
as well as:
\begin{eqnarray}
2\frac{d}{dx}\mathcal{J}_{n}^{\mu}(x)&=&\mathcal{J}_{n-1}^{\mu}(x)-\mathcal{J}_{n+1}^{\mu}(x)+\frac{2\mu}{2\mu+1}\mathcal{J}_{n+1}^{\mu+1}(x),\\
\frac{2 n}{x}\mathcal{J}_{n}^{\mu}(x)&=&\mathcal{J}_{n-1}^{\mu}(x)+\mathcal{J}_{n+1}^{\mu}(x)-\frac{2\mu}{2\mu+1}\mathcal{J}_{n+1}^{\mu+1}(x),\\
(xD_{\mu}-[n]_{\mu})\mathcal{J}_{n}^{\mu}(x)&=&-x\mathcal{J}_{n+1}^{\mu}(x)+\frac{2\mu }{2\mu+1}x\mathcal{J}_{n+1}^{\mu+1}(x).
\end{eqnarray}

The function $\mathcal{J}_{n}^{\mu}(x)$ can also writen in terms of the generalized hypergeometric function $\medskip _{1}F_{2}$  as:
\begin{eqnarray}
\mathcal{J}_{n}^{\mu}(x)&=&\frac{1}{n!}\left(\frac{x}{2}\right)^n\medskip  \medskip _{1}F_{2}
\left(
\begin{array}{ll}
\frac{1}{2}\\
n+1&\mu+\frac{1}{2}
\end{array}
;-\frac{x^2}{4}\right),
\end{eqnarray}
and it has the following asymptotic form, as $|x|\rightarrow \infty$ and $|\arg(x)|\le\frac{\pi}{2}-\varepsilon$ for some $\varepsilon>0$:
\begin{eqnarray}
\mathcal{J}_n^{\mu}(x)&\approx&\frac{\Gamma(\mu+\frac{1}{2})}{\pi}\left(\frac{x}{2}\right)^{-\mu-\frac{1}{2}}\cos\left(x-\frac{\pi}{2}(\frac{1}{2}+n+\mu)\right)\nonumber\\
&+&\frac{\Gamma(\mu+\frac{1}{2})}{\Gamma(n+\frac{1}{2})\Gamma(\frac{1}{2})}\left(\frac{x}{2}\right)^{-1}.
\end{eqnarray}
\section{Other possibilities of deformation}
\subsection{Three other versions of the deformed Bessel functions}
We give three more deformations of the classical Bessel functions, with partial analogues for the set of formulae \eqref{bessel1} to \eqref{bessel4}. A complete analogue of all formulae together does not seem to be available.\\

\noindent
\textbf{First deformation:}
\begin{equation}
J_{n}^{(1,\mu)}(x)=\sum_{k \geq 0}\frac{(-1)^k }{k![k+n]_{\mu}!}\left(\frac{x}{2}\right)^{2k+n}.
\end{equation}

\text{Generating function}:
\begin{equation}
E_{\mu}\left({\frac{xt}{2}}\right)e^{-\frac{x}{2t}}=J_0^{(1,\mu)}(x)+\sum_{n=1}^{+\infty}J_n^{(1,\mu)}(x)\big(t^n+(-1)^nt^{-n}\big).
\end{equation}

\text{Recursive relations}:
\begin{eqnarray}
xD_{\mu}J_{n}^{(1,\mu)}(x)&=&-xJ_{n+1}^{(1,\mu)}(x)+[n]_{\mu}J_{n}^{(1,\mu)}(x),
\\
2\frac{d}{dx}J_{n}^{(1,\mu)}(x)&=&\frac{1}{2\mu+1}J_{n-1}^{(1,\mu)}(x)-J_{n+1}^{(1,\mu)}(x)+\frac{\mu}{(2\mu+1)^2}x J_{n-2}^{(1,\mu+1)}(x),\\
\frac{2 n}{x}J_{n}^{(1,\mu)}(x)&=&\frac{1}{2\mu+1}J_{n-1}^{(1,\mu)}(x)+J_{n+1}^{(1,\mu)}(x)+\frac{\mu}{(2\mu+1)^2}x J_{n-2}^{(1,\mu+1)}(x)
,\\
&&\hskip -40mm (xD_{\mu}-[-n]_{\mu})J_{n}^{(1,\mu)}(x)=\frac{1}{2\mu+1}xJ_{n-1}^{(1,\mu)}(x)+\frac{\mu}{(2\mu+1)^2}x^2 J_{n-2}^{(1,\mu+1)}(x).
\end{eqnarray}

\text{Connection formula}:
\begin{equation}
 J_{n}^{(1,\mu)}(x)=\sum_{k\geq 0}\frac{(-1)^k(\mu)_k}{k!(2\mu+1)_k}x^k{J}_{n-k}(x).
\end{equation}

\noindent
\textbf{Second deformation:}
\begin{equation}
J_{n}^{(2,\mu)}(x)=\sum_{k \geq 0}\frac{(-1)^k }{[k]_{\mu}!(k+n)!}\left(\frac{x}{2}\right)^{2k+n}.
\end{equation}

\text{Generating function}:
\begin{equation}
   e^{\frac{xt}{2}}E_{\mu}\left(-\frac{x}{2t}\right)=J_0^{(2,\mu)}(x)+\sum_{n=1}^{+\infty}J_n^{(2,\mu)}(x)\big(t^n+(-1)^nt^{-n}\big).
\end{equation}

\text{Recursive relations}:
\begin{eqnarray}
 xD_{\mu}J_{n}^{(2,\mu)}(x)&=&xJ_{n-1}^{(2,\mu)}(x)+[-n]_{\mu}J_{n}^{(2,\mu)}(x),
\\
2\frac{d}{dx}J_{n}^{(2,\mu)}(x)&=&J_{n-1}^{(2,\mu)}(x)-\frac{1}{2\mu+1}J_{n+1}^{(2,\mu)}(x)+\frac{\mu}{(2\mu+1)^2}x J_{n+2}^{(2,\mu+1)}(x),\\
\frac{2 n}{x}J_{n}^{(2,\mu)}(x)&=&J_{n-1}^{(2,\mu)}(x)+\frac{1}{2\mu+1}J_{n+1}^{(2,\mu)}(x)-\frac{\mu}{(2\mu+1)^2}x J_{n+2}^{(2,\mu+1)}(x),\\
(xD_{\mu}-[n]_{\mu})J_{n}^{(2,\mu)}(x)&=&-\frac{1}{2\mu+1}xJ_{n+1}^{(2,\mu)}(x)+\frac{\mu}{(2\mu+1)^2}x^2 J_{n+2}^{(2,\mu+1)}(x).
\end{eqnarray}

\text{Connection formula}:
\begin{equation}
  J_{n}^{(2,\mu)}(x)=\sum_{k\geq 0}\frac{(\mu)_k}{k!(2\mu+1)_k}x^k{J}_{n+k}(x).
\end{equation}
\noindent
\textbf{Third deformation:}
\begin{equation}
J_{n}^{(3,\mu)}\left(x\right)=\sum_{k \geq 0}\frac{(-1)^k }{[k]_{\mu}![k+n]_{\mu}!}\left(\frac{x}{2}\right)^{2k+n}.
\end{equation}

\text{Generating function}:
\begin{equation}
E_{\mu}\left({\frac{xt}{2}}\right)E_{\mu}\left({-\frac{x}{2t}}\right)=\sum_{n=-\infty}^{+\infty}J_n^{(3,\mu)}(x)t^n.
\end{equation}
\section{Appendix: hypergeometric functions}
The generalized hypergeometric function $_pF_q$ is defined by the series 
\begin{equation}
_pF_q\left(
\begin{array}{llll}
a_1&a_2&...&a_p\\
b_1&b_2&...&b_q
\end{array};x
\right)
=\sum_{n=0}^{\infty}\frac{(a_1)_n(a_2)_n...(a_p)_n}{(b_1)_n(b_2)_n...(b_q)_n}\frac{x^n}{n!}
\end{equation}
It can be shown that the series converges for all $x$ if $p\leq q$, converges for $|x|<1$ if $p=q+1$, and diverges for all $x\neq 0$ if $p>q+1$. It is a solution of the differential equation:
\begin{equation}
\left(x\frac{d}{dx}+a_1\right)...\left(x\frac{d}{dx}+a_p\right)y-\frac{d}{dx}\left(x\frac{d}{dx}+b_1-1\right)...\left(x\frac{d}{dx}+b_q-1\right)y=0.
\end{equation}

The following differential recursive equations hold:

\begin{equation*}
\begin{array}{lll}
\left(x\frac{d}{dx}+a_1\right)\,_pF_q\left(
\begin{array}{llll}
a_1&a_2&...&a_p\\
b_1&b_2&...&b_q
\end{array};x
\right)&=&a_1\,_pF_q\left(
\begin{array}{llll}
a_1+1&a_2&...&a_p\\
b_1&b_2&...&b_q
\end{array};x
\right),\\
\left(x\frac{d}{dx}+b_1-1\right)\,_pF_q\left(
\begin{array}{llll}
a_1&a_2&...&a_p\\
b_1&b_2&...&b_q
\end{array};x
\right)&=&(b_1-1)\,_pF_q\left(
\begin{array}{llll}
a_1&a_2&...&a_p\\
b_1-1&b_2&...&b_q
\end{array};x
\right),\\
\frac{d}{dx}\,_pF_q\left(
\begin{array}{llll}
a_1&a_2&...&a_p\\
b_1&b_2&...&b_q
\end{array};x
\right)&=&\frac{a_1...a_p}{b_1...b_q}\,_pF_q\left(
\begin{array}{llll}
a_1+1&a_2+1&...&a_p+1\\
b_1+1&b_2+1&...&b_q+1
\end{array};x
\right).
\end{array}
\end{equation*}
Special cases of hypergeometric functions are for example:
\begin{equation}
e^x=\,_0F_0\left(\begin{array}{l}
-\\
-
\end{array};x
\right),
\end{equation}
\begin{equation}
(1-x)^{-a}=\,_1F_0\left(\begin{array}{l}
a\\
-
\end{array};x
\right),
\end{equation}

The Bessel function of order $\nu$ can be expressed in two different ways as a hypergeometric function:
\begin{eqnarray}
J_{\nu}(x)&=&\frac{(x/2)^{\nu}}{\Gamma(\nu+1)}\,_0F_1\left(\begin{array}{l}
-\\
\nu+1
\end{array};-\frac{x^2}{4}
\right)\\
&=&\frac{e^{-ix}(x/2)^{\nu}}{\Gamma(\nu+1)}\,_1F_1\left(\begin{array}{l}
\nu+1/2\\
2\nu+1
\end{array};2ix
\right).
\end{eqnarray}
For a detailed account, see for example \cite{Erdelyi, lebedev}, or any textbook on special functions.

 \end{document}